\documentstyle[12pt]{article}
\newtheorem{prop}{Proposition}
\newtheorem{th}{Teorem}

\newtheorem{cor}{Corollary}
\begin{document}
\author{A.A.Davydov}
\title{Quasitriangular structures on cocommutative Hopf algebras}
\maketitle
\begin{abstract}
The article is devoted to the describtion of quasitriangular structures 
(universal R-matrices) on cocommutative Hopf algebras. It is known that such 
structures are concentrated on finite dimensional Hopf subalgebras. 

In particular, quasitriangular structure on group algebra is defined by the 
pairs of normal inclusions of an finite abelian group and by invariant 
bimultiplicative form on it. The structure is triangular in the case of 
coinciding inclusions and skewsymmetric form.

The nonstandart $\lambda$-structure on the representation ring of finite 
group, corresponding to the triangular structure on group ring, is described.
\end{abstract}
\tableofcontents
\section{Introduction}

The notion of quasitriangular Hopf algebra was introduced by V.G.Drinfeld
in connection with the so-called quantum Yang-Baxter equation, which 
appears in different contexts of mathematical physics.

Quasitriangular Hopf algebra is a pair, which consists of a Hopf algebra and
invertible element of its tensor square (universal R-matrix), satisfying some
additional conditions. 

From categorical point of view quasitriangular Hopf algebras are characterized 
by the fact that its representations categories (categories of modules having 
finite dimensional over ground field) are rigid and quasitensor.

The present article is devoted to the describtion of quasitriangular structures 
(universal R-matrices) on cocommutative Hopf algebras over algebraically closed 
field of characteristic zero. It is used that any quasitriangular structure 
are concentrated on some finite demensional Hopf subalgebra \cite{ra}. 
Since finite dimensional cocommutative Hopf algebras over algebraically closed 
field of characteristic zero are group algebras, the special attention is 
devoted to the case of this class of Hopf algebras.

It is shown that quasitriangular structures on group algebra is defined by the 
pairs of normal inclusions of an finite abelian group and by invariant 
bimultiplicative form on it (theorem \ref{gra}). The triangular structures 
corresponds to the case of coinciding inclusions and skewsymmetric form.

It is well known that the tensor structure on monoidal category allows to 
define external powers of any object and gives a $\lambda$-structure on its 
Grothendieck ring. In the last section 
the nonstandart $\lambda$-structures on the representation ring of finite 
group, corresponding to the triangular structures on group ring, is described 
(theorem \ref{ada}). It is appears that these $\lambda$-structure depends only 
of central involution of the group, which is defined by skewsymmetric 
bimultiplicative form.

The work was partially supported by Russian Fund of Fundamental Research, 
grant n 96-01-00149.

\section{Tensor categories}

Let $k$ be a ground field. All categories and functors will be $k$-linear. 

The {\em monoidal category} is a category {\cal G} \ with a bifunctor
\begin{displaymath}
\otimes :{\cal G} \times {\cal G} \longrightarrow {\cal G} \qquad (X,Y) \mapsto 
X \otimes Y
\end{displaymath}
which called {\em tensor (or monoidal) product}. This functor must be equiped 
with a functorial collection of isomorphisms (so-called {\em associativity 
constraint})
\begin{displaymath}
\varphi_{X,Y,Z} : X \otimes (Y \otimes Z) \rightarrow (X \otimes Y) \otimes Z 
\qquad \mbox{for any} \quad X,Y,Z \in {\cal G}
\end{displaymath}
which satisfies to the following {\em pentagon axiom}:
$$(X\otimes\varphi_{Y,Z,W})\varphi_{X,Y\otimes Z,W}(\varphi_{X,Y,Z}\otimes W) = 
\varphi_{X,Y,Z\otimes W}\varphi_{X\otimes Y,Z,W}.$$

Consider two tensor products of objects $X_{1},...,X_{n}$ from ${\cal G}$ with 
an arbitrary arrangement of the brackets.
The coherence theorem \cite{mcl:cwm} states that the pentagon axiom implies the 
existence of a unique isomorphism between them, which is the composition of the 
associativity constraints. This fact allows us to omit brackets in the tensor 
products. 
\newline
An object 1 together with the functorial isomorphisms
$$\rho_{X}:X \otimes 1 \rightarrow X \quad
\lambda_{X}:1 \otimes X \rightarrow X$$
in a monoidal category ${\cal G}$ \ is called a {\em unit} if $\lambda_1 = 
\rho_1$ and 
$$\lambda_{X\otimes Y} = \lambda_{X}\otimes I:1\otimes X\otimes Y\to 
X\otimes Y,\quad
\rho_X\otimes I = I\otimes\lambda_Y:X\otimes 1\otimes Y\to X\otimes Y,$$
$$\rho_{X\otimes Y} = I\otimes\rho_{Y}:X\otimes Y\otimes 1\to X\otimes Y$$
for any $X,Y\in\cal G$. 

It is easy to see that the unit object is unique up to isomorphism.
We will suooise additionly that the endomorphisms ring $End_{\cal G}(1)$ of the 
unite object coincides with the ground field. 

A {\em monoidal functor} between monoidal categories ${\cal G}$
and ${\cal H}$ is a functor $F : {\cal G} \longrightarrow
{\cal H}$ , which is equipped with the functorial collection
of isomorphisms (the so-called {\em monoidal structure})
\begin{displaymath}
F_{X,Y} : F(X \otimes Y) \rightarrow F(X) \otimes F(Y)
\qquad \mbox{for any} \quad X,Y \in {\cal G}
\end{displaymath}
for which 
$$F_{X,Y\otimes Z}(I\otimes F_{Y,Z}) = F_{X\otimes Y,Z}(F_{X,Y}\otimes I)$$
for any objects $X,Y,Z \in {\cal G}$.

A morphism $f : F \rightarrow G$ of monoidal functors $F$ and $G$ is called 
{\em monoidal} if 
$$G_{X,Y}f_{X \otimes Y} = (f_{X} \otimes f_{Y})F_{X,Y}$$
for any $X,Y \in {\cal G}$.

The tensor product allows to correspond to any functor $F:\cal G\to H$ from 
monoidal category $\cal G$ the collection of functors from the cartesian powers 
of the category $\cal G$
$$F^{\otimes n}:{\cal G}^{\times n}\to{\cal H}\quad  
F^{\otimes n+1}(X_{1},...,X_{n+1}) = F(...(X_{1} \otimes...) \otimes X_{n+1})$$
The collection of the endomorphisms algebras of tensor powers of a monoidal 
functor F can be equipped with the structure of cosimplicial complex 
$E(F)_{*} = End(F^{\otimes *})$.
\newline
The image of the coface map
$${\partial}_{n+1}^{i} : End(F^{\otimes n}) \rightarrow End(F^{\otimes n+1}) 
\qquad i = 0,...,n+1$$
of the endomorphism $\alpha \in End(F^{\otimes n})$ has the following 
specialization on the objects $X_{1},...,X_{n+1}$:
$${\partial}_{n+1}^{i}(\alpha )_{X_{1},...,X_{n+1}} = \left\{
\begin{array}{l}
{\phi}_{0}(I_{X_{1}} \otimes \alpha_{X_{2},...,X_{n+1}}){\phi}^{-1}_{0}, \quad 
i = 0 \\
{\phi}_{i}(\alpha_{X_{1},...,X_{i} \otimes X_{i+1},...,X_{n+1}}){\phi}^{-1}_{i}, 
\quad 1 \leq i \leq n \\
\alpha_{X_{1},...,X_{n}} \otimes I_{X_{n+1}}, \quad i = n+1
\end{array}
\right.,$$
here ${\phi}_{i}$ is the unique isomorphism between $F(...(X_{1} \otimes...) 
\otimes X_{n+1}) = F^{\otimes n+1}(X_{1},...,X_{n+1})$ and $F(...(X_{1} 
\otimes...) \otimes (X_{i} \otimes X_{i+1})) \otimes ...) \otimes X_{n+1})$.

The specialization of the image of the coboundary map
$${\sigma}_{n-1}^{i} : End(F^{\otimes n}) \rightarrow End(F^{\otimes n+1}) 
\qquad i = 0,...,n-1$$
is
$${\sigma}_{n-1}^{i}(\alpha)_{X_{1},...,X_{n-1}} = 
{\alpha}_{X_{1},...,X_{i},1,X_{i+1},...,X_{n-1}}.$$

Let us note that the monoidality of the automorphism $f\in Aut(F)$ of monoidal 
functor $F$ can be presented in the following form: 
$${\partial}^{1}_{2}(f) = {\partial}^{0}_{2}(f){\partial}^{2}_{2}(f).$$

Two monoidal functors will be called by {\em twisted forms} of each other if
they are isomorphic as functors. It this casae the differ only by the monoidal 
structure. The ratio of its monoidal structures is an automorphism 
$\alpha\in Aut(F^{\otimes 2})$ of the tensor square of one of these functors. 
The fact that the composition $F_{X,Y}\alpha_{X,Y}$ is a monoidal structure on 
the functor $F$ is equivalent to the condition: 
\begin{equation}\label{tmf}
{\partial}^{0}_{3}(\alpha ){\partial}^{2}_{3}(\alpha ) =  
{\partial}^{3}_{3}(\alpha ){\partial}^{1}_{3}(\alpha ).
\end{equation}
So twisted forms of monoidal functors $F$ correspond to the automorphisms 
$\alpha\in Aut(F^{\otimes 2})$, which satisfies (\ref{tmf}). 
We will call this automorphisms by {\em 2-cocycles} of the monoidal functor $F$ 
and by $Z^2 (F)$ we will denote the set of 2-cocyles of $F$. 

The twisted form of the funtor $F$ corresponding to 2-cocycle 
$\alpha\in Z^2 (F)$ will be denoted by $F(\alpha )$.  

Monoidal category $\cal G$ with the unite object is {\em rigid} if it is 
equipped by the contravariant equivalence of categories  
$(\ )^* :\cal G\to\cal G$ ({\em dualization}) together with the functorial 
collection of morphisms 
$$\kappa_X :1\to X\otimes X^*,\quad ev_X :X^*\otimes X\to 1,$$ 
for which the compositions 
$$X\to^{\kappa\otimes I}X\otimes X^*\otimes X\to^{I\otimes ev}X,$$
$$X^*\to^{I\otimes\kappa}X^*\otimes X\otimes X^*\to^{ev\otimes I}X^*$$
are identical. 

It is easy to see that for any monoidal functor $F$ between rigid categories 
the composition $F((\ )^* )$ is canonically isomorphic to $F(\ )^*$. 
The isomorphism can be presented as the composition 
$$F(X^* )\to^{I\otimes\kappa_{F(X)}}F(X^* )\otimes F(X)\otimes F(X)^*\to^{F_{
X^* ,X}\otimes I}F(X^*\otimes X)\otimes F(X)^*\to^{F(ev_{X})\otimes I}F(X)^*.$$

The monoidal category $\cal G$ is {\em quasitensor} if it is eqiupped by the 
functorial collection of isomorphisms ({\em commutativity constraint})
$c_{X,Y}:X\otimes Y\to Y\otimes X$, for which 
$$c_{X,Y\otimes Z} = (Y\otimes c_{X,Z})(c_{X,Y}\otimes Z)\quad
c_{X\otimes Y,Z} = (c_{X,Z}\otimes Y)(X\otimes c_{Y,Z}).$$ 
Quasitensor category is {\em tensor} if $c_{Y,X}c_{X,Y} = I$ for any 
$X,Y\in\cal G$.

Monoidal functor $F:\cal G\longrightarrow\cal H$ between (quasi)tensor 
categories is {\em (quasi)tensor} if 
$$F_{Y,X}F(c_{X,Y}) = c_{F(X),F(Y)}F_{X,Y},$$
for any $X,Y\in\cal G$.

For the rigid quasitensor category the square of the dualization is isomorphic 
to the identity functor $\nu :Id\to (\ )^{**}$. The isomorphism can be presented 
as the functorial composition 
$$X\to^{I\otimes\kappa_{X^*}}X\otimes X^*\otimes X^{**}\to^{c\otimes I}
X^*\otimes X\otimes X^{**}\to^{ev\otimes I}X^{**}.$$
The morphism $\nu$ allows to define the automorphism $\gamma$ of identity 
functor $Id_{\cal G}$: 
$$\gamma_X :X\to^{\nu_X}X^{**}\to^{\nu_{X^{**}}}X^{****}\to^{(\nu_{X^*})^*}
X^{**}\to^{\nu^{-1}_{X}}X.$$

Let us remark that the square of the dualization $(\ )^{**}$ is a monoidal 
functor. The morphism $\nu$ is monoidal in the case of tensor category. In the 
general case it defines 2-cocycle $\sigma\in Z^2 (id_{\cal G})$ of the identity 
functor 
$$\sigma_{X,Y}:X\otimes Y\to^{\nu_X\otimes\nu_Y}X^{**}\otimes Y^{**}\to 
(X\otimes Y)^{**}\to^{\nu_{X\otimes Y}^{-1}}X\otimes Y.$$
The proof of the following proposition can be extracted from \cite{re}.
\begin{prop}
1. $\sigma_{X,Y} = c_{Y,X}c_{X,Y},$ 

in particular $\sigma=1$ for the tensor category. 

2. $\sigma^2 = \partial (\gamma).$
\end{prop}

The isomorphism $\nu$ allows to define the trace of any endomorphism of any 
object. Let $f\in End_{\cal G}(X)$ is an endomorphism of the object $X$ of 
rigid quasitensor category $\cal G$. The {\em trace} of the endomorphism $f$ is 
an element of the ground field $Tr_{X}(f)\in k$, for which the composition 
$$1\to^{\kappa_{X}} X\otimes X^*\to^{\nu_X f\otimes I}X^{**}\otimes X^*
\to^{ev_{X^*}}1$$
coincides with $Tr_{X}(f)Id_1$.

The trace is additive by the definition:
$$Tr_{X}(f+g) = Tr_{X}(f)+Tr_{X}(g),\quad\forall f,g\in End_{\cal G}(X).$$
In the case of tensor category the trace is also multiplicative: 
$$Tr_{X\otimes Y}(f\otimes g) = Tr_{X}(f)Tr_{Y}(g),\quad\forall 
f\in End_{\cal G}(X), g\in End_{\cal G}(Y).$$

The {\em rank} of the object $X$ is a trace of its identity morphism 
$rk(X) = Tr_{X}(I_{X})$.  

Using the morphism $\nu$ we can define an automorphism 
$\mu (F) = \nu_{\cal H}^{-1}F(\nu_{\cal G})$ of any monoidal functor 
$F:\cal G\longrightarrow\cal H$ between quasitensor categories which will be 
called {\em Markov automorphism}. Let us note that $\mu (F) = 1$, iff $F$ is 
a quasitensor functor.  
\begin{prop}\label{trf}
Let $F:\cal G\longrightarrow\cal H$ be monoidal functor between quasitensor 
categories and $f\in End_{\cal G}(X)$. Then 
1. $F(\sigma_{\cal G}) = \partial (\mu (F))\sigma_{\cal H},$ 

in particular $F(\sigma ) = \partial (\mu (F)),$ in the case of tensor 
category $\cal H$

2. $Tr_{X}(f) = Tr_{F(X)}(\mu (F)_{X}F(f)).$
\end{prop}
Proof:

The first equation can be obtained by replacing $F(\nu_X )$ on 
$\nu_{F(X)}\mu(F)_X$ in the expression of $F(\sigma_{X,Y})$. 

As was mentined above the object $F(X^* )$ is canonically isomorphic to 
$F(X)^*$ and this isomorphism identifies $F(\kappa_X )$ with $\kappa_{F(X)}$ and 
$F(ev_X )$ with $ev_{F(X)}$. Hence the composition 
$$1\to^{F(\kappa_{X})} F(X\otimes X^* )\to^{F(\nu_X f\otimes I)}F(X^{**}
\otimes X^*)\to^{F(ev_{X^*})}1,$$
wich is equal to $F(Tr_{X}(f)) = Tr_{X}(f)$, is isomorphic to the following 
$$1\to^{\kappa_{F(X)}} F(X)\otimes F(X^* )\to^{\nu_{F(X)}F(f)\otimes I}
F(X^{**})\otimes F(X^*)\to^{ev_{F(X)^*}}1.$$ 
Now it is enough to note that $F(\nu_X ) = \nu_{F(X)}\mu(F)_X$. $\Box$

\section{Quasitriangular Hopf algebras}

This section contains the definitions and general properties of quasitriangular 
structures on the Hopf algebras. 

The {\em Hopf algebra} is an associative unitary algebra $H$ with the algebra
homomorphisms:
$$\Delta :H\otimes H\rightarrow H\quad(\mbox{coproduct}),$$
$$\varepsilon :H\rightarrow k\quad(\mbox{counite})$$ 
and with the antihomomorphism of algebras:
$$S :H\rightarrow H\quad(\mbox{antipode}),$$
which satisfies to the following conditions: 
$$ (I\otimes\Delta )\Delta = (\Delta\otimes I)\Delta\quad
(\mbox{coassociativity}),$$
$$ (I\otimes\varepsilon )\Delta = (\varepsilon\otimes I)\Delta = I\quad
(\mbox{counitarity}),$$
$$ (I\otimes S)\Delta = (S\otimes I)\Delta = \varepsilon.$$
Here $I$ is an identical map. 

An element $g\in H$ for which $\Delta (g) = g\otimes g$ will be called 
{\em grouplike}. It is easy to verify that $g\quad \varepsilon (g) = 1$ and 
$S(g) = g^{-1}$ for grouplike $g$. The set of grouplike elements of Hopf 
algebra $H$ form a group $G(H)$ under multiplication. 

Coproduct allows to define the structure of $H$-module on the tensor product 
$M\otimes_{k}N$ of two $H$-modules: 
$$h*(m\otimes n) = \Delta(h)(m\otimes n)\qquad h\in H, m\in M, n\in N.$$ 
Coassociativity of coproduct implices that the standart associativity constrint 
of vector spaces 
$$\varphi :L\otimes (M\otimes N)\rightarrow (L\otimes M)\otimes N \qquad 
\varphi (l\otimes (m\otimes n) = (l\otimes m)\otimes n$$
is $H$-linear. 

The counite defines $H$-module structure on the ground field $k$ and the 
counite axiom means that this is a unite object. 

Finally the antipode of Hopf algebra allows to define (left) $H$-module 
structure on the dual space $M^* = Hom(M,k)$ of any (left) $H$-module $M$: 
$$h*l(x) = l(S(h)x),\qquad\forall h,x\in H, l\in M^* .$$
It follows from the antipode axiom that for $H$-module $M$ finite dimensional 
over ground field the Kazimir inclusion $k\to M\otimes M^*$ and the evaluation 
map $M^*\otimes M\to k$ are $H$-linear. 

By another words the category $Rep(H)$ of {\em representations} (that are (left) 
$H$-modules which are finite dimensional over ground field) of Hopf algebra $H$ 
is rigid monoidal. It is easy to see that for any {\em homomorphism of Hopf 
algebras} $f:H\to F$ (algebra homomorphism for which 
$\Delta f = (f\otimes f)\Delta$) the restriction functor 
$f^*:Rep(F)\to Rep(H)$ is monoidal with trivial monoidal structure. For example, 
the forgetful functor $Rep(H)\to k-mod$, which coincides with the restriction 
functor of the inclusion $k\to H$, is monoidal. 

Accordingly to the previous section nontrivial monoidal structures on the 
restriction functor correspond to some invertible elements of its emdomorphisms 
complex. 

The {\em cobar complex} of Hopf algebra $H$ is a tensor algebra 
$\oplus H^{\otimes n}$ with coface maps given by 
$${\partial}^{i}_{n}:H^{\otimes n-1}\longrightarrow H^{\otimes n},$$
$${\partial}^{i}_{n}(h_{1}\otimes ...\otimes h_{n}) =
\left\{
\begin{array}{cl}
1\otimes h_{1}\otimes ...\otimes h_{n},& i=0\\
h_{1}\otimes ...\otimes\Delta (h_{i})\otimes ...\otimes h_{n},& 1\leq i\leq n\\
h_{1}\otimes ...\otimes h_{n}\otimes 1,& i=n+1
\end{array}
\right.,$$
and codegeneration maps
$${\sigma}^{i}_{n}(h_{1}\otimes ...\otimes h_{n+1}) = h_{1}\otimes ...\otimes
\varepsilon (h_{i})\otimes ...\otimes h_{n+1})).$$
The following proposition was proved in \cite{dav}. 
\begin{prop}
Endomorphisms complex of the restriction functor 
$f^{*}:Rep(F)\rightarrow Rep(H)$ coincides with the subcomplex 
$C_{F^{\otimes *}}(\Delta (f(H))$ of centralizers of the image of coproduct in 
the cobar complex of Hopf algebra $F$. 

An isomorphism is given by the map 
$$H^{\otimes n}\longrightarrow End(F^{\otimes n}),$$
sending the element $x\in H^{\otimes n}$ to the endomorphism of multiplication 
by $x$.
\end{prop}

In particular twisted forms of the restriction functor $f^*:Rep(F)\to Rep(H)$ 
correspond to the invertible elements $F\in C_{F^{\otimes 2}}(\Delta (f(H))$, 
satisfying to the condition: 
\begin{equation}\label{cc}
(1\otimes F)(I\otimes\Delta )(F) = (F\otimes 1)(\Delta\otimes I)(F).
\end{equation}

Hopf algebra is {\em cocommutative} if $t\Delta = \Delta$, where 
$t:H\otimes H\rightarrow H\otimes H$ denotes the permutation of tensor factors. 
Let us note that the antipode of cocommutative Hopf algebra is involutive: 
$S^2 = I$.

It is not hard to verify that finite dimensional cocommutative Hopf algebra 
over the algebraically closed field of characteristic zero is a group algebra of
its grouplike elements \cite{lar}. 

The following notion can be regarded as a generalization of the cocommutativity.
The pair $(H,R)$, which consists of the Hopf algebra $H$ and the invertible 
element $R\in H^{\otimes 2}$, is {\em quasitriangular} Hopf algebra if the 
following conditions are fulfilled:
\begin{equation}\label{inv}
Rt\Delta (h) = \Delta (h)R \quad \mbox{for any}\quad h\in H,
\end{equation}
\begin{equation}\label{fe}
(I\otimes \Delta)(R) = R_{13}R_{12},
\end{equation}
\begin{equation}\label{se}
(\Delta\otimes I)(R) = R_{13}R_{23}.
\end{equation}
The element $R$ will be called an {\em universal $R$-matrix}.

quasitriangular Hopf algebra will be called {\em triangular} if 
$$RR_{21} = 1 \quad\mbox{(unitarity of the universal $R$-matrix)}.$$ 

The proof of the following proposition is contained in the Drinfeld's paper 
\cite{dr}.
\begin{prop}\label{prm}
The universal $R$-matrix $R$ in the quasitriangular Hopf algebra satisfies to 
the following conditions:
$$(\varepsilon\otimes I)(R) = (I\otimes\varepsilon )(R) = 1,$$
$$(S\otimes I)(R) = (I\otimes S)(R) = R^{-1},\quad (S\otimes S)(R) = R,$$
$$R_{12}R_{13}R_{23} = R_{23}R_{13}R_{12}\quad
\mbox{(quantum Yang-Baxter equation)}.$$
\end{prop}

The universal $R$-matrix allows to define quasitensor structure on the 
representations category $Rep(H)$:
$$c_{M,N}:M\otimes N\to N\otimes M,\quad m\otimes n\mapsto R(n\otimes m).$$
Let us note that the unitarity of the $R$-matrix is equivalent to the tensority 
of the corresponding structure.

Let $(H,R)$ and $(H',R')$ be two quasitriangular Hopf algebras and  
$f:H\to H'$ is a homomorphism of Hopf algebras. It is easy to see that the 
twisted form of the restriction functor $f^*$, corresponding to 2-cocycle 
$F\in Z^2 (f^* )$, is a quasitensor functor iff $R'F = t(F)R$. 

The element $u = (S\otimes I)(R_{21})$ will be called {\em Markov element} of 
the universal $R$-matrix $R$. Let us note that Markov authomorphism of 
forgetful functor $Rep(H)\to k-mod$ coincides with the multiplication by 
$u^{-1}$. 

The proof of the following proposition, which borrows from the Drinfeld's paper
\cite{dr}, can be also concluded from the results of first article.
\begin{prop}
Markov element $u$ satisfies to the following conditions:

1) $u$ is invertible, and $u^{-1} = (S^{-1}\otimes S)(R_{21})$,

2) $S^2 = Ad(u)$, that is $S^2 (h) = uhu^{-1}$ for any $h\in H$,
in particular, $u$ is central for cocommutative $H$,
 
3) $\Delta (u) = (R_{21}R)^{-1}(u\otimes u) = (u\otimes u)(R_{21}R)^{-1}$,
in particular, $u$ is grouplike for unitary $R$-matrix $R$,

4) $\Delta (g) = (g\otimes g)$ and $g = uS(u)^{-1}$.
\end{prop}

\section{Minimal quasitriangular Hopf algebras}

This section contains the describtion of minimal quasitriangular Hopf algebras
which belongs to Radford \cite{ra}.

Any bivector $R\in H^{\otimes 2}$ provides two finite dimensional subspaces of 
$H$
$$H_l = \{ (I\otimes l)(R),\quad l\in H^* \}\qquad
H_r = \{ (l\otimes I)(R),\quad l\in H^* \}.$$
Here $H^*$ denotes the space $Hom(H,k)$ of linear functions on $H$.
\newline
Let us remark that $R$ lies in $H_l\otimes H_r$ by the consctruction.
\newline
Moreover, the bivector $R$ defines two bijective maps: 
$$\alpha :{H_r}^*\rightarrow H_l\qquad \alpha (l) = (I\otimes l)(R),$$
$$\beta :{H_l}^*\rightarrow H_r\qquad \beta (l) = (l\otimes I)(R),$$
for which
$$\alpha^* = \beta.$$
Indeed, by the definition of dual map $m(\alpha^* (l)) = l(\alpha(m))$ for any 
$l,m\in H^*$. Hemce
$$m(\alpha^* (l)) = l(\alpha(m)) = (l\otimes m)(R) = m(\beta (l)),
\quad\mbox{for any}\quad l,m\in H^* .$$

It is appear that the spaces $H_l$ and $H_R$, defined by an universal 
$R$-matrix, are Hopf subalgebras in $H$.

Here we will need the notions of dual Hopf algebra and adjoint action.

It is easy to see that the dual maps to the product, coproduct 
and antipode are coproduct, product and antipode on the dual space $H^*$ for
finite dimensional Hopf algebra $H$:
$$\Delta_H^* = \mu_{H^*},\quad\mu_H^* = \Delta_{H^*},\quad i_H^* = 
\epsilon_{H^*},\quad S_H^* = S_{H^*}.$$
This Hopf algebra will be called {\em dual}. 

Let us define the {\em adjoint} action of Hopf algebra $H$ on itself by setting
$x^h = \sum_{(h)}h_{(0)}xS(h_{(1)})$ for $x\in H$. This action induce also the 
action of $H$ on $H^*$:
$$l^h (x) = l(x^{S(h)}),\qquad\forall l\in H^* , h,x\in H.$$
We will call subspace in $H$ {\em normal} if it is invariant respectively the
adjoint action.
\begin{th}\label{sup}
For any universal $R$-matrix on the Hopf algebra $H$

subspaces $H_l,H_r$ are finite dimensional Hopf subalgebras of $H$,

the map $\alpha$ is an antihomomorphism of algebras and homomorphism of 
coalgebras. 

In the case of unitary $R$-matrix subalgebras 
$H_l$ and $H_r$ coincides and the map 
$\alpha$ satisfies to the condition 
$$\alpha_{R}^* = S\alpha_{R}.$$

If $H$ is cocommutative Hopf algebra, then the maps $\alpha$ and $\beta$ are
$H$-invariant:
$$\alpha (l^h ) = \alpha (l)^h ,\qquad 
\beta (l^h ) = \beta (l)^h ,\quad\forall h\in H, l\in H^* .$$
In particular, $H_l, H_r$ are normal subalgebras. 
\end{th}
Proof:

Let $l$ and $m$ be two linear functions on $H$.
Applying $I\otimes l\otimes m$ to both sides of the equation (\ref{fe}), we 
will receive that the product of any two elements from $H_l$ lies in $H_l$: 
$$(I\otimes l)(R)(I\otimes m)(R) = (I\otimes l\otimes m)(R_{13}R_{12}) =
(I\otimes l\otimes m)(I\otimes \Delta )(R).$$
The equality $(I\otimes\epsilon )(R) = 1$ (proposition \ref{prm}) means that 
the identity also lies in $H_l$. So $H_l$ is a subalgebra of $H$.

To verify that $H_l$ is a subcoalgebra of $H$ it is enough to apply
$I\otimes l$ to both sides of the equation (\ref{se}):
$$(\Delta\otimes l)(R) = (I\otimes I\otimes l)(\Delta\otimes I)(R) = 
(I\otimes l)(R_{13}R_{23}).$$

The fact that $H_r$ is a Hopf subalgebra can be proved analogously.

Antihomomorphity of the map $\alpha$ also follows from the equation (\ref{fe}). 
Let $l$ and $m$ be two linear functions from $H_r^*$. Then
$$\alpha(l*m) = (I\otimes l*m)(R) = 
(I\otimes l\otimes m)(I\otimes \Delta )(R) =$$
$$(I\otimes l\otimes m)(R_{13}R_{12}) = (I\otimes m)(R)(I\otimes l)(R) = 
\alpha(m)\alpha(l).$$ 
The equality (\ref{se}) implices that the map $\alpha$ is a coalgebra 
homomorphism: 
$$\Delta (\alpha_{R}(l)) = (\Delta\otimes l)(R) = 
(I\otimes I\otimes l)(\Delta\otimes I)(R) = (I\otimes l)(R_{13}R_{23}).$$
$$(I\otimes \Delta (l))t_{23}(R\otimes R) = 
(\alpha_{R}\otimes\alpha_{R})(\Delta (l)),$$
where $t_{23}$ is a permutation of second and third tensor fuctors. 

The unitarity condition $Rt(R) = 1$ and the equation $(I\otimes S)(R) = R^{-1}$ 
(proposition \ref{prm}) implice that $H_l$ coincides with $H_r$:
$$(I\otimes l)(R) = (l\otimes I)(t(R)) = (l\otimes I)(R^{-1}) = 
(lS\otimes I)(R).$$ 

Let us check the equality $\alpha^* = S\alpha$:
$$S\alpha (l) = (S\otimes l)(R) = (I\otimes l)(R^{-1}) = $$
$$(I\otimes l)(t(R)) = \alpha^* (l).$$

In the case of cocommutative Hopf algebra $H$ the condition (\ref{inv}) is 
equivalent to the $H$-invariance of universal $R$-matrix:
$$R\Delta (h) = \Delta (h)R \quad \mbox{for any}\quad h\in H,$$
which is equivalent to the condition:
$$\sum_{(h)}(h_{(0)}\otimes 1)R(S(h_{(1)})\otimes 1) = 
\sum_{(h)}(1\otimes S(h_{(0)}))R(1\otimes h_{(1)}).$$ 
Applying $(1\otimes l)$ to the previous equality, we will obtain 
$H$-invariance of the map $\alpha$:
$$\alpha (l^h ) = \alpha (l)^h.$$
$H$-invariance of $\beta$ can be verified analogously. 
$\Box$

So we can define by the universal $R$-matrix $R$ in the Hopf algebra $H$ the 
finite dimensional Hopf algebra $F (=H_l )$ and two Hopf algebra inclusions 
$i:F\to H, j:F^*\to H$. Moreover $R$ coincides with the image $(i\otimes j)(C)$ 
of Kazimir element $C\in F\otimes F^*$. To characterize the image of the map
$i\otimes j$ it is convinient to use the notion of quantum double.

{\em Quantum double} $D(A)$ of (finite dimensional) Hopf algebra $A$ \cite{dr} 
is unique quasitriangular Hopf algebra satisfying to the following conditions:

1) $D(A)$ contains $A$ and $A^{*op}$ as subalgebras. Here $A^{*op}$ denotes the 
dual Hopf algebra with opposite multiplication.

2) Universal $R$-matrix $R$ is the image of Kazimir element $C\in A\otimes A^*$ 
under natural inclusion of vector spaces 
$A\otimes A^*\rightarrow D(A)\otimes D(A)$.

3) The multiplication $A\otimes A^{*op}\rightarrow D(A)$ is a bijective. 

The quasitriangular Hopf algebra $(H,R)$ will be called {\em minimal} if its 
not contains proper quasitriangular subalgebras.

By means of quantum double the theorem \ref{sup} can be stated in the following 
form. 
\begin{cor}
For any minimal quasitriangular Hopf algebra $(H,R)$ there are finite 
dimensional Hopf algebra $F$ and the surjective homomorphism of quasitriangular 
algebras $(D(F),R)\to (H,R)$, which restrictions on $F$ and $F^*$ are injective. 
\end{cor}

Let us cosider the case of triangular Hopf algebra in more detail. 
We will call universal $R$-matrix $R$ {\em nondegenerated} if the corresponding 
map $\alpha_{R}:H^*\to H$ is a bijective. 
\begin{th}
For any unitary universal $R$-matrix $R$ in the Hopf algebra $H$ there is a Hopf 
subalgebra $F$, such that $R\in F^{\otimes 2}$ and $R$ is nondegenerated on $F$.
\end{th}
Proof:

Let us denote by $F$ the image of the map $\alpha_{R}$. By the construction  
the bivector $R$ lies in $F^{\otimes 2}$. 

Let us show that $\alpha_{R}$ decomposes into a product of the maps
$$H^*\rightarrow^{i^*} F^*\rightarrow^{\alpha} F\rightarrow^{i}H,$$ 
where $i$ is a standart inclusion. 
Write $\alpha_{R}$ as the composition of standart inclusion and projection: 
$$H^*\rightarrow^{\pi} F\rightarrow^{i}H.$$
Since by the theorem \ref{sup} 
$\alpha_{R}^* = S\alpha_{R}$, the map $\pi$ differs from
$i^*$ by an automorphism of Hopf algebra $F$. 
$\Box$

For cocommutative algebra $H$ the finite dimensional Hopf algebra $F$ is 
commutative and cocommutative. It is follows from the describtion of 
cocommutative Hopf algebras over the algebraically closed fields of 
characteristic zero \cite{lar} that $F$ is a group algebra of finite abelian 
group. 
 
\section{Group algebras}
Let us remind that our ground field $k$ is algebraically closed of 
characteristic zero. 

\begin{th}\label{gra}
Universal $R$-matrix in the group algebra $k[G]$ of a group $G$ is defined by 
the pair of normal inclusions $i,j:A\to G$ of finite abelian group, which induce 
the same $G$-module structure on $A$, and by nondegenerated bimultiplicative 
$G$-invariant form $\beta :\hat A\otimes\hat A :\rightarrow k^*$. 

The universal $R$-matrix defined by this data have the following form: 
$$R = \frac{1}{{\vert A\vert}^2}\sum_{a,b\in A}\sum_{\chi ,\xi\in\hat A}
\beta (\chi ,\xi )\chi (a)\xi (b)(i(a)\otimes j(b)).$$

Unitary $R$-matrix corresponds to the case of coinciding inclusions $i=j$ and 
skewsymmetric form. 

The Markov element $u$ of the unitary $R$-matrix, defined by the form 
$\beta$, is a central involution of the group $G$ and satisfies to the equation, 
which defines it uniqually:
\begin{equation}\label{me}
\chi (u) = \beta (\chi ,\chi )\qquad\mbox{for any}\quad\chi\in\hat A.
\end{equation}
\end{th}
Remark. Here $\hat A$ denotes the characters group $Hom(A,k^* )$ of the group 
$A$. Bimultiplicative form $\beta :\hat A\otimes\hat A :\rightarrow k^*$ is a 
homomorphism from tensor product to the group of nonzero elements of the 
ground field. Nondegeneracy of $\beta$ means that the map 
$$\hat A :\rightarrow {\hat{\hat A}}\quad \chi\mapsto \beta (\chi ,?)$$
is bijective.

Proof:

Let $R$ be a $R$-matrix in the group algebra $k[G]$. 
Accordingly to the theorem \ref{sup} there are finite abelian group $A$ and two 
inclusions $i,j:A\to G$ with normal images such that the homomorphism of Hopf 
algebras $\alpha_R :k[G]^*\to k[G]$ decomposes into the product 
$$k[G]^*\to^{i^*}k[A]^*\to^{\alpha}k[A]\to^{j}k[G],$$
where $\alpha$ is an isomorphism. 

For finite abelian group $A$ any isomorphism of Hopf algebras 
$\alpha :k[A]^*\to k[A]$ is defined by the nondegenerated bimultiplicative form 
$\beta :\hat A\otimes \hat A\to k^*$. Indeed the Hopf algebra 
$k[A]^*$ coincides with the group Hopf algebra $k[\hat A]$ of the characters 
group of $A$. The isomorphism of group Hopf algebras 
$\alpha :k[\hat A]\to k[A]$ is defined by an isomorphism of groups 
$\hat A\to A$. It is enough to note that such isomorphisms are the same as 
nondegenerated bimultiplicative forms $\beta :\hat A\otimes \hat A\to k^*$. 

By the theorem \ref{sup} the inclusions $i,j:A\to G$, defined by unitary 
$R$-matrix, coincide, and isomorphism $\alpha :k[A]^*\to k[A]$ satisfies to the 
condition $\alpha^* = S\alpha$. The direct cheking shows the equivalence of this 
condition to skewsymmetricity of the corresponding form. 

By the definition of Markov element:
$$u = \mu (S\otimes I)(R) = 
\frac{1}{{\vert A\vert}^2}\sum_{a,b\in A}\sum_{\chi ,\xi\in\hat A}
\beta (\chi ,\xi )\chi (a)\xi (b)a^{-1}b =$$
$$\frac{1}{{\vert A\vert}^2}\sum_{a,c\in A}\sum_{\chi ,\xi\in\hat A}
\beta (\chi ,\xi )\chi\xi (a)\xi (c)c =
\frac{1}{\vert A\vert}\sum_{c\in A}\sum_{\chi\in\hat A}
\beta (\chi ,\chi )^{-1}\chi (c)c.$$
Now let us use the condition (\ref{me}):
$$\xi (u) = \frac{1}{\vert A\vert}\sum_{c\in A}(\sum_{\chi\in\hat A}
\beta (\chi ,\chi )^{-1}\chi\xi (c) = 
\beta (\xi ,\xi^{-1}) = \beta (\xi ,\xi).$$
The fact that the Markov element is an involution follows from 
skewsymmetricity of the form $\beta$.
$\Box$

Example.

Let $A$ be a cyclic group of order 2 with the generator $u$. There is unique 
nondegenerated skewsymmetric bimultiplicative form 
$\beta :\hat A\otimes \hat A\to k^*$. The corresponding $R$-matrix has the 
following form: 
$$R_u = \frac{1}{2}(1\otimes 1 + 1\otimes u + u\otimes 1 - u\otimes u).$$
The representations category of the group $A$ with tensor structure, defined by 
this $R$-matrix, is equivalet to the tensor category ${\cal K}{\it oz}$, whose 
objects are $\bf Z/2\bf Z$-graded vector spaces with commutativity constraint, 
defined by Kozul rule of sings \cite{dm}:
$$x\otimes y\mapsto (-1)^{\tilde x\tilde y}y\otimes x$$
for homogeneous $x,y$. Here $\tilde x$ denotes the degree of the element $x$. 
$\Box$

Let us note that for arbitrary unitary $R$-matrix $R$ in the group algebra 
$k[G]$ there is not tensor functor from the representations category of $G$ with 
tensor structure defined by the $R$ to the category of vector spaces. However 
the following fact take place.
\begin{prop}
For any tensor structure on the representatins category ${\cal R}{\it ep}_k (G)$ 
of finite group $G$ there is a tensor functor to the category 
${\cal K}{\it oz}$. 
\end{prop}
Proof:

Let $R$ be unitary $R$-matrix in $k[G]$, which defines the given tensor 
structure on the representation category ${\cal R}{\it ep}_k (G)$. Let $A$ be 
an abelian normal subgroup of $G$ and $\beta$ - nondegenerated skewsymmetric 
bimultiplicative form on the characters group $\hat A$, which corresponds to the 
$R$-matrix $R$.

Let us note that the Markov element $u$ of $R$-matrix $R$ belongs to the 
subgroup $A$. We can define bimultiplicative form $\beta_u$ on $\hat A$ as the 
composition 
$$\beta_u :\hat A\otimes\hat A\to {\hat {\langle u\rangle}}\otimes
{\hat {\langle u\rangle}}\to k^*$$
of the restriction on the cyclic subgroup $\langle u\rangle$, generated by the 
element $u$, and unique nondegenerated skewsymmetric bimultiplicative form on 
$\hat{\langle u\rangle}$.

The forms $\beta ,\beta_u$ are related as follows: 
$$\beta_u (\chi ,\chi ) = \chi (u) = \beta (\chi ,\chi )\quad
\forall\chi\in\hat A.$$
Hence its ratio lies in the kernel of the map 
$$Hom(\Lambda^2 (\hat A),k^* )\to Hom({\hat A},k^* )\simeq A,$$
which sends skewsymmetric bimultiplicative form $\beta$ to the multiplicative 
form $l(\chi ) = \beta (\chi ,\chi )$. 
Let us note that for any abelian group $B$ the map 
$Hom(\Lambda^2 (B),k^* )\to Hom(B,k^* )$ can be complemented to the exact 
sequence 
$$H^2 (B,k^* )\to Hom(\Lambda^2 (B),k^* )\to Hom(B,k^* ),$$ 
which first map sends the class of 2-cocycle $\gamma$ to the form 
$\beta (a,b) = \gamma (a,b)\gamma (b,a)^{-1}$ \cite{br,hu}.
In our case we can find the 2-cocycle $\gamma\in Z^2 (A,k^* )$ such that 
$$\beta (\chi ,\xi )\gamma (\xi ,\chi ) = 
\beta_u (\chi ,\xi )\gamma (\chi ,\xi ).$$
Let us define an element $F_{\gamma}\in k[G]^{\otimes 2}$ by 
$$F_{\gamma} = \sum_{a,b\in A}(\sum_{\chi ,\xi\in\hat A}
\gamma (\chi ,\xi )\chi (a)\xi (b)a\otimes b.$$ 
By the construction it satisfies to the conditions: 
$$Rt(F_{\gamma}) = R_{u}F_{\gamma},$$ 
$$(u\otimes u)F_{\gamma} = F_{\gamma}(u\otimes u),$$
$$(1\otimes F_{\gamma})(I\otimes\Delta )(F_{\gamma}) = 
(F_{\gamma}\otimes 1)(\Delta\otimes I)(F_{\gamma}).$$
Hence it defines the tensor structure of the restriction functor 
$$res^{G}_{\langle u\rangle}:{\cal R}{\it ep}_k (G)\to {\cal K}{\it oz}.$$
$\Box$

\section{Nonstandart $\lambda$-structures on the character rings of finite 
groups}
Let us remind the definition of the $\lambda$-ring. 

Let $A$ be a commutative ring without additive torsion. 
The set of maps $\lambda^i :A\rightarrow A\qquad i\geq 0$, satisfying to the 
conditions
\begin{equation}\label{dlo}
\lambda^0 (x)=1\quad \lambda^1 (x)=x\quad \lambda^i (x+y) = 
\sum _{s+t=i} \lambda^s (x)\lambda^t (y)\qquad\mbox{for any}\quad x,y\in A,
\end{equation}
will be called {\em $\lambda$-operations} or {\em external powers}. 
It is convinient to write the conditions (\ref{dlo}) by means of generating 
series $\lambda_t (x)=\sum_{i\geq 0}\lambda^i (x)$
$$\lambda_t (x+y) = \lambda_t (x)\lambda_t (y)\quad\mbox{для любых}\quad
x,y\in A.$$ 
The generating series is also usefull to define another types of operations:

{\em $\sigma$-operations} or ({\em symmetric powers}) 
$\sigma^i :A\rightarrow A\qquad i\geq 0$, with generating series 
$$\sigma_t (x) = \lambda_{-t}(x)^{-1}$$ 
and {\em Adams operations} $\psi^i :A\rightarrow A\qquad i\geq 1$, 
with generating series 
$$\psi_{-t}(x) = -t\frac{{\lambda '}_t (x)}{\lambda_t (x)}.$$ 
The generating series of these operations satisfies to the conditions:
$$\sigma_t (x+y) = \sigma_t (x)\sigma_t (y)\quad
\psi_t (x+y) = \psi_t (x) + \psi_t (y)
\quad\mbox{for any}\quad x,y\in A.$$ 

The ring $A$ with the collection of $\lambda$-operations will be called 
{\em $\lambda$-ring}, if its Adams operations are ring homomorphisms. 

The tensor structure on the category of complex representations
${\cal R}{\it ep}_k(G)$ of finite group $G$, which corresponds to the unitary 
$R$-matrix $R$, allows to define the structure of the $\lambda$-ring on the 
characters ring $R_k(G)$ (the Grothendieck ring of the category 
${\cal R}{\it ep}_k (G)$). 

The action of the symmetric group $S_n$ on the tensor power $X^{\otimes n}$ of 
representation $X$, defined by the tensor structure, provides the decomposition 
of $X^{\otimes n}$ into the direct sum of isotypical components, labeled by the 
irreducible representations of $S_n$. The {\em external power} $\Lambda^n (X)$ 
is an isotypical component $Hom_{S_n}(sign,X^{\otimes n})$ of the nontrivial one 
dimensional representation $sign$. 

The external power $\lambda^{n}_{R}[X]$ of the class $[X]$ of the 
representation $X$ is by definition the class of its external power 
$[\Lambda^{n}_{R}(X)]$. Using the additivity of external powers these operations 
can be extended to representations ring. 

The $\lambda$-structure of representations ring can be expressed in more simple 
form by means of characters.

We can identify the representations ring with the subring (ring of characters) 
in the ring of complex functions of conjugate classes by corresponding to the 
representation the traces of elements of the group. Adams operations have the 
following form on the characters \cite{at}:
$$\psi^i (x)(g) = x(g^i ) \quad\mbox{for any}\quad x\in R_{\bf C}(G), 
g\in G.$$

We will use the following proposition in the proof of the general statement of 
this article. 
\begin{prop}\label{adams}
Let $\cal G$ be an abelian $k$-linear tensor category over an algebraically 
closed fields $k$ of characteristic zero. Then for any monoidal automorphism $g$ 
of the identity functor of the category $\cal G$
$$Tr(g_{\psi^k (X)}) = Tr(g^{k}_{X})\quad\forall X\in R_k (G).$$
\end{prop}
Proof:

Using the property $\psi^{nm} = \psi^n\psi^m$ of the Adams operations 
\cite{gr,kn} we can suppose that $k=p$ is a prime number. 

Let us use the representation of Adams operation $\psi^p$ as a difference of 
cyclic operations $\psi^p = c^{p}_{1} - c^{p}_{\varepsilon}$, where 
$\varepsilon$ is a nontrivial degree $p$ root of unity \cite{ker,ati}.

By the definition, the value of the cyclic operation, corresponding to degree 
$p$ root of unity $\varepsilon$, on the class of object $X$ is a class of the 
image of the projector on $p$-th tensor power $X^{\otimes p}$ 
$$c^{p}_{\varepsilon}(X) = \frac{1}{p}\sum_{i=0}^{p-1}\varepsilon^i\tau^i 
)(X^{\otimes p}),$$ 
where $\tau$ is a cycle of length $p$ in the symmetric group $S_p$. 

By the additivity of trace we have
$$Tr(g_{c^{p}_{\varepsilon}(X)}) = 
\frac{1}{p}\sum_{i=0}^{p-1}\varepsilon^i Tr(\tau^i g_{X}^{\otimes p}),$$ 
and using lemma 7.2 from \cite{de} we can write
$Tr_{X^{\otimes p}}(\tau^i g_{X}^{\otimes p}) = Tr_{X}(g^{p}_{X})$. 
Hence the trace
$$Tr(g_{c^{p}_{\varepsilon}(X)}) = Tr_{X}(g^{p}_{X}) 
\frac{1}{p}\sum_{i=0}^{p-1}\varepsilon^i$$
equals zero, for nontrivial root of unity $\varepsilon$ and 1, if
$\varepsilon = 1$.
$\Box$

\begin{th}\label{ada}
The $\lambda$-structure of the representation ring $R_k (G)$ of finite group 
$G$, corresponding to the unitary $R$-matrix $R\in k[G]^{\otimes 2}$, depends 
only of its Markov element $u\in G$. 

The Adams operations of this $\lambda$-structure have the following form:
\begin{equation}\label{nls}
\psi^{k}_{u}(\chi )(g) = \chi (u^{k+1}g^k ),\quad\forall \chi\in R_k (G), 
g\in G.
\end{equation}
\end{th}
Proof:

Additivity of Adams operations allows to verify the equation (\ref{nls}) only 
for the characters of representations. 
Let $\chi$ be the character of the representation $X\in {\cal R}{\it ep}(G)$. 
By the definition $\chi (g)$ coincides with the trace $Tr_{\omega (X)}(g_{X})$, 
where $\omega :{\cal R}{\it ep}(G)\to k-mod$ is forgetful functor.
By the proposition \ref{trf} we can replace $Tr_{\omega (X)}(g_{X})$ by 
$Tr_{X}(u_{X}g_{X})$. In particular 
$$\psi^{k}_{u}(\chi )(g) = Tr_{\omega (X)}(g_{\psi^{k}_{u}(X)}) = 
Tr_{\psi^{k}_{u}(X)}((ug)_{\psi^{k}_{u}(X)}).$$
By the proposition \ref{adams} the previous expansion coinsides with the 
following 
$$Tr_{X}((ug)^{k}_{X}) = Tr_{\omega (X)}((u^{k+1}g^{k})_{X}) = 
\chi (u^{k+1}g^k ).$$
$\Box$

\end{document}